\documentclass[preprint,12pt]{elsarticle}

\usepackage{geometry}
\usepackage{amsmath,amssymb,theorem}
\usepackage{float}
\usepackage{hyperref}
\usepackage{graphicx}
\usepackage{longtable}
\usepackage{caption}
\newtheorem{thm}{Theorem}[section]
\newtheorem{claim}{Claim}[section]
\newtheorem{conj}[thm]{Conjecture}
\newtheorem{lemma}[thm]{Lemma}
\newtheorem{coro}[thm]{Corollary}
\newtheorem{prob}[thm]{Problem}
\newtheorem{case}{Case}
\newtheorem{subcase}{Case}[case]
\newtheorem{subsubcase}{Case}[subcase]

\journal{acta mathematicae applicatae sinica}

\begin{document}

\begin{frontmatter}

\title{Hadwiger's conjecture for some graphs with independence number two}

\author[a,b]{Tong Li}
\ead{litong@amss.ac.cn}

\author[a,b]{Qiang Zhou}
\ead{zhouqiang2021@amss.ac.cn}

\address[a]{Academy of Mathematics and Systems Science, Chinese Academy of Sciences, Beijing, 100190, China}

\address[b]{University of Chinese Academy of Sciences, Beijing, 100049, China}

\begin{abstract}

Let $h(G)$ denote the largest $t$ such that $G$ contains $K_t$ as a minor and $\chi(G)$ be the chromatic number of $G$ respectively. In 1943, Hadwiger conjectured that $h(G) \geq \chi(G)$ for any graph $G$. In this paper, we prove that Hadwiger's conjecture holds for \emph{H-free} graphs with independence number two, where $H$ is one of some specified graphs.

\end{abstract}

\begin{keyword}
Hadwiger's conjecture \sep graph minor \sep independence number.


\end{keyword}

\end{frontmatter}

\section{Introduction}
All graphs considered in this paper are finite, undirected and simple. Let $V(G)$, $|G|$, $E(G)$, $e(G)$, $\overline{G}$, $\delta(G)$, $\Delta(G)$, $\omega(G)$, $\alpha(G)$ and $\chi(G)$ denote the vertex set, number of vertices, edge set, number of edges, complement, minimum degree, maximum degree, clique number, independence number, and chromatic number of a graph $G$, respectively. Let $e(G,H)$ denote the number of edges between two graphs $G$ and $H$. A vertex set $A$ is \emph{complete} to a vertex set $B$, if every vertex in $A$ is adjacent to every vertex in $B$. For $X \subseteq V(G)$, Let $G[X]$ denote the subgraph of $G$ induced by $X$, $G - X$ denote the subgraph $G[V(G) \backslash X]$ or simply $G - y$ when $X = \{ y \}$, respectively. Let $N_{G}(X)$ denote the set of vertices adjacent to $X$ in $G - X$. A graph $H$ is called a \emph{minor} of $G$ if $H$ is obtained from $G$ by some sequence of vertex deletions, edge deletions, and edge contractions. Let $h(G)$ denote the largest $t$ such that $G$ contains $K_t$ as a minor.

In 1943, Hadwiger \cite{hadwiger1943klassifikation} conjectured that:

\begin{conj}[\cite{hadwiger1943klassifikation}]\label{HC}
For every graph $G$, $h(G) \geq \chi(G)$.
\end{conj}

Hadwiger's conjecture is one of the most famous open problems in graph theory. The Four-color Theorem is a special case when we restrict $h(G) \leq 4$. Hadwiger \cite{hadwiger1943klassifikation} and Dirac \cite{dirac1952property} proved the conjecture holds for $h(G) \leq 3$ independently. Robertson, Seymour and Thomas \cite{robertson1993hadwiger} proved this conjecture is equivalent to the Four-color Theorem when $h(G) \leq 5$. So far, Hadwiger's conjecture is still open when $h(G) \geq 6$.

An equivalent statement of Hadwiger's conjecture is that every graph with no $K_{t+1}$ minor is $t$-colorable. Table~\ref{table1} lists some previous studies:

\begin{longtable}[c]{|c|c|l|}
\hline
Author(s)           & Year & \multicolumn{1}{c|}{Results}                                                                                                   \\ \hline

Delcourt, Postle \cite{delcourt2021reducing}  & 2021 & \begin{tabular}[c]{@{}l@{}}Every graph with no $K_{t+1}$ minor is\\ $O(t \log \log t)$-colorable.\end{tabular}                 \\ \hline
Wagner \cite{wagner1960bemerkungen}          & 1960 & \begin{tabular}[c]{@{}l@{}}Every graph with no $K_{5}^{-}$ (remove one\\ edge from $K_{5}$) minor is 4-colorable.\end{tabular} \\ \hline
Dirac \cite{dirac1964structure}           & 1964 & \begin{tabular}[c]{@{}l@{}}Every graph with no $K_{6}^{-}$ minor is\\ 5-colorable.\end{tabular}                                \\ \hline
Jakobsen \cite{jakobsen1971weakenings}\cite{jakobsen1970homomorphism} &
  \begin{tabular}[c]{@{}c@{}}1970\\ 1971\end{tabular} &
  \begin{tabular}[c]{@{}l@{}}Every graph with no $K_{7}^{-}$ minor is\\ 7-colorable;\\ Every graph with no $\mathcal{K}_{7}^{-2}$ (remove two\\ edges from $K_{7}$) minor is 6-colorable.\end{tabular} \\ \hline
Rolek, Song \cite{rolek2017coloring} &
  2017 &
  \begin{tabular}[c]{@{}l@{}}Every graph with no $K_{8}^{-}$ minor is\\ 9-colorable;\\ Every graph with no $\mathcal{K}_{8}^{-2}$ minor is\\ 8-colorable.\end{tabular} \\ \hline
Lafferty, Song \cite{lafferty2022every} &
  2022 &
  \begin{tabular}[c]{@{}l@{}}Every graph with no $\mathcal{K}_{8}^{-4}$ minor is\\ 7-colorable;\\ Every graph with no $\mathcal{K}_{9}^{-6}$ minor is\\ 8-colorable.\end{tabular} \\ \hline
Rolek \cite{rolek2020graphs}          & 2020 & \begin{tabular}[c]{@{}l@{}}Every graph with no $\mathcal{K}_{9}^{-2}$ minor is\\ 10-colorable.\end{tabular}                    \\ \hline
Kostochka \cite{kostochka2014minors} &
  2014 &
  \begin{tabular}[c]{@{}l@{}}Every graph with no $K_{s,t}$ minor is\\ $(s+t-1)$-colorable when\\ $t > 5(200s \log_{2}(200s))^{3}$.\end{tabular} \\ \hline
Kostochka, Prince \cite{kostochka2010dense} & 2010 & \begin{tabular}[c]{@{}l@{}}Every graph with no $K_{3,t}$ minor is\\ $(t+2)$-colorable when $t \geq 6300$.\end{tabular}         \\ \hline

\caption{}
\label{table1}
\end{longtable}

Seymour \cite{seymour2016hadwiger} pointed out that graphs $G$ with $\alpha(G) = 2$ are particularly interesting. He believed that Hadwiger's conjecture is probably true if it holds for graphs with independence number two.

\begin{prob}\label{HC2}
Does $h(G) \geq \chi(G)$ hold for every graph $G$ with $\alpha(G) \leq 2$ \rm ?
\end{prob}

One idea to study Problem~\ref{HC2} is to find a big complete minor, according to the following theorem:

\begin{thm}[\cite{plummer2003special}]\label{HC2iff}
For every graph $G$ with $\alpha(G) \leq 2$, $h(G) \geq \chi(G)$ if and only if $h(G) \geq \lceil |G| / 2 \rceil$.
\end{thm}

Recently, Norin and Seymour \cite{norin2022dense} proved that for any graph $G$ with $\alpha(G) \leq 2$, there exists a minor with $\lceil |G| / 2 \rceil$ vertices and with fewer than $1/76$ of all possible edges missing.

\begin{figure}[H]
\centering
\includegraphics[height=9cm,width=13.5cm]{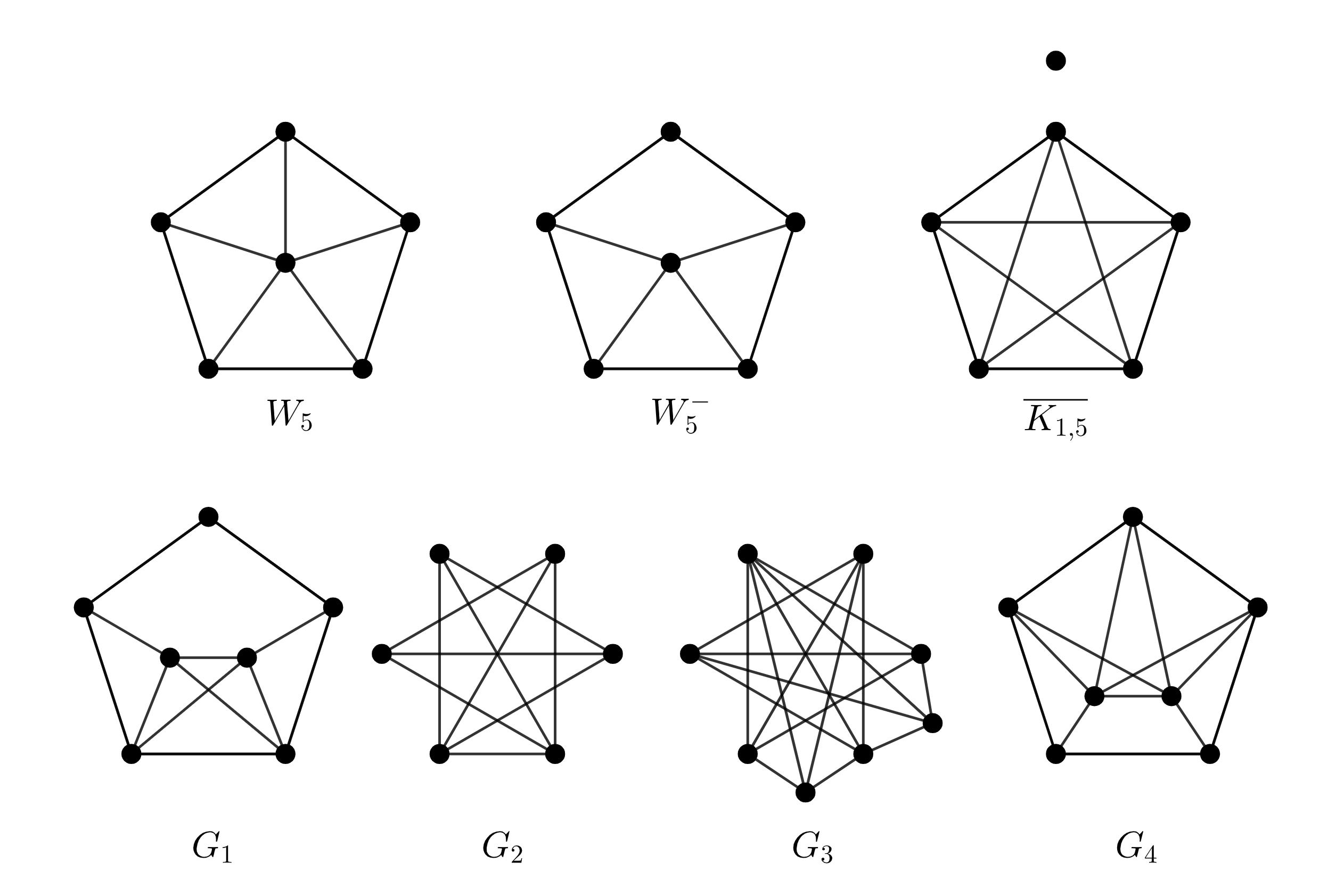}
\caption{Graph $H$ in previous studies}
\label{previous}
\end{figure}

We can also consider Problem~\ref{HC2} for $H$-free graphs, that is, graphs without $H$ as an induced subgraph. Table~\ref{table2} and Figure~\ref{previous} list previous studies of Problem~\ref{HC2} for $H$-free graphs:

\begin{longtable}[c]{|c|c|c|}
\hline
Author(s)           & Year & \multicolumn{1}{c|}{Results}                                                                                                   \\ \hline

Plummer, Stiebitz, Toft \cite{plummer2003special} & 2003 & \begin{tabular}[c]{@{}c@{}}$|H| \leq 4$, or $H = C_{5}$,\\ or $H = G_{1}$ in Figure~\ref{previous}.\end{tabular}                                   \\ \hline
Kriesell \cite{kriesell2010seymour}               & 2010 & \begin{tabular}[c]{@{}c@{}}$|H| \leq 5$,\\ or $H = W_{5}^{-}$ in Figure~\ref{previous}.\end{tabular}                                                   \\ \hline
Bosse \cite{bosse2019note}                        & 2019 & \begin{tabular}[c]{@{}c@{}}$H = W_{5}$ in Figure~\ref{previous},\\ or $H = \overline{K_{1,5}}$ in Figure~\ref{previous},\\ or $H = K_{7}$.\end{tabular} \\ \hline
Carter \cite{carter2022hadwiger}           & 2022+ & \begin{tabular}[c]{@{}c@{}}$H$ is any of the 33 graphs in \cite{carter2022hadwiger}\\ or $H = K_{8}$. (computer-assisted)\end{tabular}                \\ \hline
Carter \cite{carter2022hadwiger}           & 2022+ & \begin{tabular}[c]{@{}c@{}}$H = G_{2}$, $G_{3}$ or $G_{4}$ in Figure~\ref{previous}.\end{tabular}                              \\ \hline
\caption{}
\label{table2}
\end{longtable}

\begin{figure}[H]
\centering
\includegraphics[height=6cm,width=9.5cm]{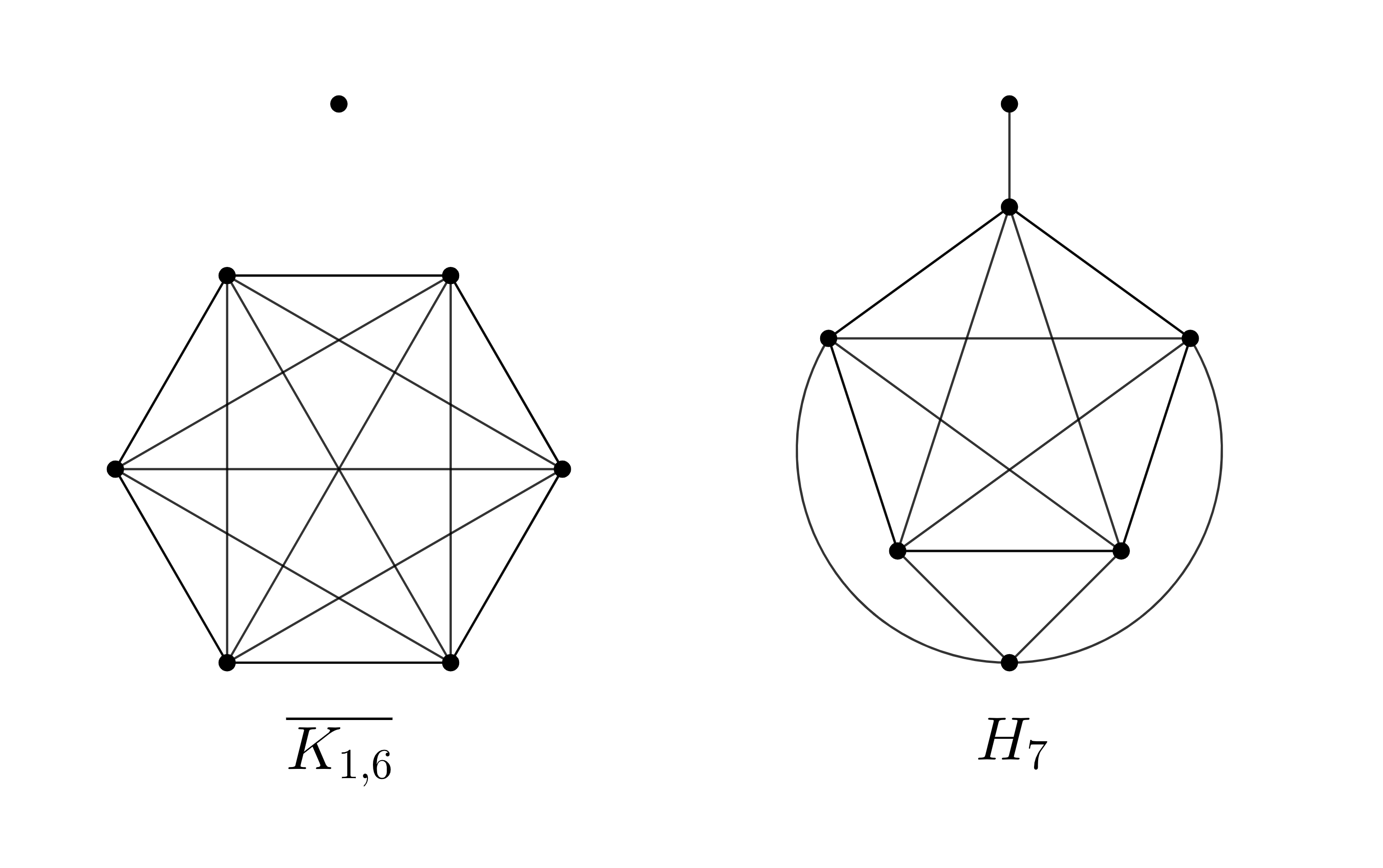}
\caption{Graph $H$ in this paper}
\label{new}
\end{figure}

In particular, Carter \cite{carter2022hadwiger} used some algorithms to verify Problem~\ref{HC2} for $H$-free graphs when $H$ is $K_{8}$ or any of the 33 graphs in the paper, and proved for $H=G_{2},G_{3},G_{4}$ mathematically. In this paper, we prove Problem~\ref{HC2} hold for more $H$ mathematically, generalize the results of Plummer, Stiebitz, Toft, Kriesell, Bosse and Carter:

\begin{thm}\label{T1}
Let $H = \overline{K_{1,6}}$ in Figure~\ref{new}. For every $H$-free graph $G$ with $\alpha(G) \leq 2$, $h(G) \geq \chi(G)$.
\end{thm}

Since $\overline{K_{1,6}}$-free is equivalent to $\delta(G) \geq |G|-6$, we can also get:

\begin{coro}
For every graph $G$ with $\alpha(G) \leq 2$, if $\delta(G) \geq |G|-6$, then $h(G) \geq \chi(G)$.
\end{coro}

\begin{thm}\label{T2}
Let $H = H_{7}$ in Figure~\ref{new}. For every $H$-free graph $G$ with $\alpha(G) \leq 2$, $h(G) \geq \chi(G)$.
\end{thm}

\begin{thm}\label{T3}
For every $K_{8}$-free graph $G$ with $\alpha(G) \leq 2$, if $\delta(G) \neq 19$, then $h(G) \geq \chi(G)$.
\end{thm}

Here we limit $\delta(G) \neq 19$ as we do not know how to
 handle the case when $|G|=27$, $\delta(G)=19$, $\Delta(G)=20 \ \text{or} \ 21$.

\section{Preliminaries}
The following theorems and lemmas are used in our proofs:

\begin{thm}[\cite{bosse2019note}]\label{W5}
For every $K_{7}$-free graph $G$ with $\alpha(G) \leq 2$, $h(G) \geq \chi(G)$.
\end{thm}

\begin{thm}[\cite{chudnovsky2012packing}]\label{packingseagull}
Let $G$ be a graph with $\alpha(G) \leq 2$. If
\[
\omega(G) \geq
\begin{cases}
   |G|/4, & \text{if $|G|$ is even}\\
   (|G|+3)/4, & \text{if $|G|$ is odd}
\end{cases},
\]
then $h(G) \geq \chi(G)$.\\
\end{thm}

\begin{lemma}\label{lemma1}
Let $G$ be a minimal counterexample to Problem~\ref{HC2} with respect to $|V(G)| + |E(G)|$, then $|G| = 27$ or $|G| \geq 29$.
\end{lemma}

\begin{pf}
By Theorem~\ref{W5}, $G$ must contain a $K_{7}$, then by Theorem~\ref{packingseagull},
\[
7 \leq \omega(G) <
\begin{cases}
   |G|/4, & \text{if $|G|$ is even}\\
   (|G|+3)/4, & \text{if $|G|$ is odd}
\end{cases},
\]
we have $|G| = 27$ or $|G| \geq 29$.\qed
\end{pf}

An edge $e$ is a \emph{dominating edge} of graph $G$ if $e$ is adjacent to all vertices in $G \backslash e$. A matching $M$ is \emph{connected} if its edges are adjacent to each other. A matching $M$ is a \emph{dominating matching} if each edge in $M$ is adjacent to all vertices in $G \backslash M$.(A dominating edge is a special connected dominating matching.)

\begin{lemma}\label{lemma2}
Let $G$ be a minimal counterexample to Problem~\ref{HC2} with respect to $|V(G)| + |E(G)|$, then $G$ does not contain a connected dominating matching.
\end{lemma}

\begin{pf}
If $G$ contains a connected dominating matching $M$, then $G \backslash M$ has a complete minor of $\lceil |G \backslash M| / 2 \rceil$ vertices. Together with the contractions of every matching edges in $M$, we get a complete minor of $\lceil |G| / 2 \rceil$ vertices from $G$, a contradiction.\qed
\end{pf}\\

\section{Proof of Theorem~\ref{T1}}
Let $G$ be a minimal counterexample to Theorem~\ref{T1} with respect to $|V(G)| + |E(G)|$. That is to say, $G$ is $\overline{K_{1,6}}$-free with $\alpha(G) \leq 2$, but $h(G) < \chi(G)$.

Since $G$ is $\overline{K_{1,6}}$-free, then we can choose a vertex $v$ with $d_{G}(v) \geq |G| - 6$. Let $A = V(G)\backslash N_{G}[v]$, then $|A| \leq 5$. For any vertex $u \in N_{G}(v)$, $u$ is not adjacent to at least one vertex in $A$, otherwise $e_{uv}$ is a dominating edge. By Lemma~\ref{lemma1}, $|G| \geq 27$, it follows that there are at least 21 non-edges between $A$ and $N_{G}(v)$. By pigeonhole principle, there exists a vertex $w \in A$ is not adjacent to at least 5 vertices in $N_{G}(v)$. Since $\alpha(G)=2$, the five vertices form a clique. Then $w$, $v$ and these 5 vertices form an induced $\overline{K_{1,6}}$, a contradiction.\\

\section{Proof of Theorem~\ref{T2}}
\begin{figure}[H]
\centering
\includegraphics[height=5cm,width=13cm]{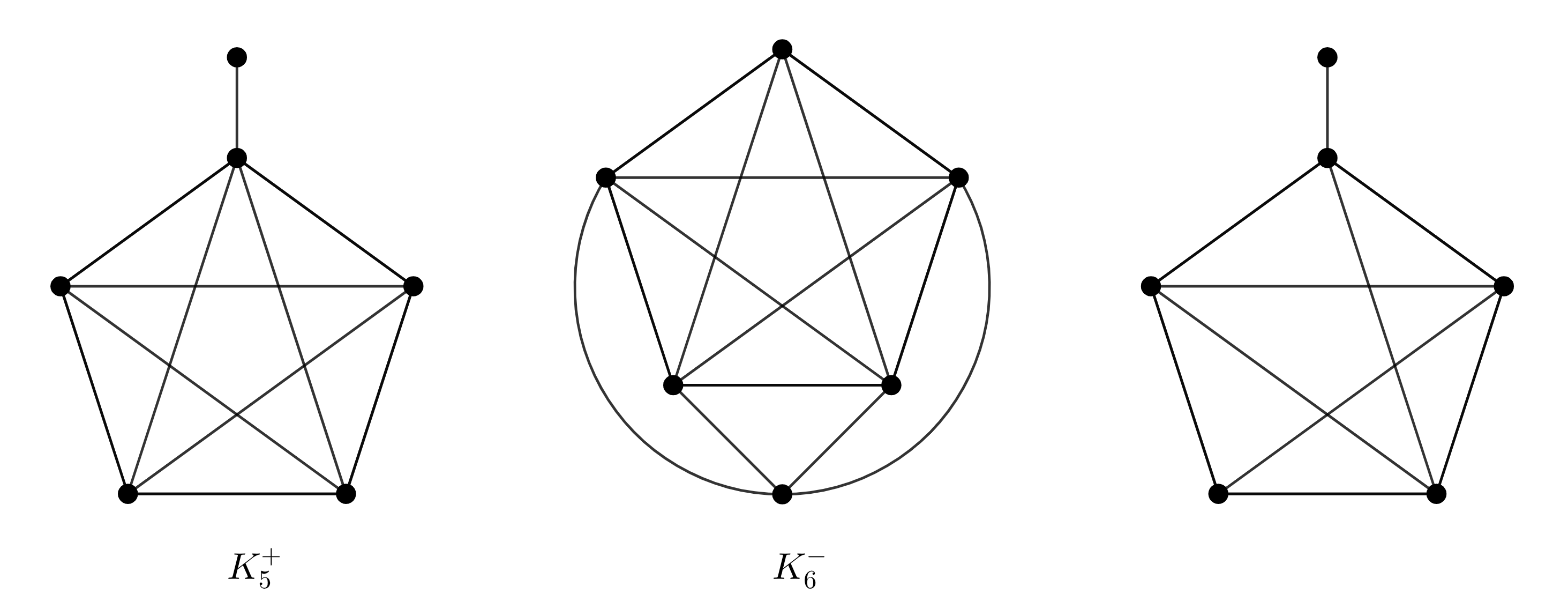}
\caption{Some induced subgraphs of $H_{7}$}
\label{induced_subgraph}
\end{figure}

Figure~\ref{induced_subgraph} lists some induced subgraphs of $H_{7}$, which are not induced subgraphs of any graph in Figure~\ref{previous}.

We first prove that Theorem~\ref{T2} holds for $H = K_{5}^{+}$. In fact, the other three graphs are simple corollaries of this case.

Let $G$ be a minimal counterexample to Theorem~\ref{T2} with respect to $|V(G)| + |E(G)|$. That is to say, $G$ is $H$-free with $\alpha(G) \leq 2$, but $h(G) < \chi(G)$. By Lemma~\ref{lemma1} we have $|G| = 27$ or $|G| \geq 29$.

By Theorem~\ref{T1}, $G$ contains an induced $\overline{K_{1,6}}$(actually $\overline{K_{1,5}}$ is enough to prove Theorem~\ref{T2}). So there exists a minimal vertex set $T$ such that $G - T$ has exactly two components $F_{1}$ and $F_{2}$, where $F_{1}$ contains a $K_{6}$ (so $|F_{1}| \geq 6$) and $|F_{2}| \geq 1$. $F_{1}$ and $F_{2}$ must be cliques since $\alpha(G) \leq 2$. For any vertex $v_{T} \in T$, neither $N_{F_{1}}(v_{T})$ nor $N_{F_{2}}(v_{T})$ is empty due to the minimality of $T$.

\begin{claim}\label{claim1}
For any vertex $v_{T} \in T$, $2 \leq |N_{F_{1}}(v_{T})| \leq 3$.
\end{claim}

\begin{pf}
If $|N_{F_{1}}(v_{T})| = 1$ for some vertex $v_{T} \in T$, then $v_{T}, N_{F_{1}}(v_{T})$ and any other four vertices in $F_{1}$ form an induced copy of $H$, a contradiction.

If $|N_{F_{1}}(v_{T})| \geq 4$ for some vertex $v_{T} \in T$, then $v_{T}$, any four vertices in $N_{F_{1}}(v_{T})$ and any one vertex in $N_{F_{2}}(v_{T})$ form an induced copy of $H$, a contradiction.
\end{pf}

\begin{claim}\label{claim2}
$T$ is complete to $F_{2}$.
\end{claim}
\begin{pf}
By Claim~\ref{claim1} and $|F_{1}| \geq 6$, any vertex in $T$ is not complete to $F_{1}$, so $T$ is complete to $F_{2}$ since $\alpha(G) \leq 2$.
\end{pf}

\begin{claim}\label{claim3}
$1 \leq |F_{2}| \leq 3$.
\end{claim}
\begin{pf}
If $|F_{2}| \geq 4$, then any one vertex $v_{T} \in T$, any four vertices in $F_{2}$ and any one vertex in $N_{F_{1}}(v_{T})$ form an induced copy of $H$, a contradiction.
\end{pf}

\begin{claim}\label{claim4}
$|F_{1}| = 6$.
\end{claim}
\begin{pf}
If $|F_{1}| \geq 7$, then any two vertices in $T$ have at most six common neighbors in $F_{1}$ by Claim~\ref{claim3}. In other words, there exists one vertex in $F_{1}$ which is not adjacent to either of these two vertices. So these two vertices are adjacent to each other since $\alpha(G) \leq 2$. By the arbitrary choice of these two vertices, we can get that $T$ is a clique.

Now $F_{1}$ and $T \cup F_{2}$ are two disjoint cliques in $G$ and $G = F_{1} \cup T \cup F_{2}$. So there exists a clique of at least $\lceil |G| / 2 \rceil$ vertices. By Theorem~\ref{HC2iff}, $h(G) \geq \chi(G)$, a contradiction.\\
\end{pf}

Next we will divide our proof into three cases:

\begin{case}
$|F_{2}| = 3$.
\end{case}

We first prove that for any two vertices $v_{T}$ and $v_{T}'$ in $T$, $N_{F_{1}}(v_{T}) = N_{F_{1}}(v_{T}')$ if and only if $v_{T}$ is adjacent to $v_{T}'$.

On the one hand, if $v_{T}$ is adjacent to $v_{T}'$ but $N_{F_{1}}(v_{T}) \neq N_{F_{1}}(v_{T}')$, then we can choose a vertex $v \in N_{F_{1}}(v_{T}) \triangle N_{F_{1}}(v_{T}') = N_{F_{1}}(v_{T}) \cup N_{F_{1}}(v_{T}') - N_{F_{1}}(v_{T}) \cap N_{F_{1}}(v_{T}')$. Thus, $\{ v \} \cup \{ v_{T} \} \cup \{ v_{T}' \} \cup F_{2}$ can form an induced copy of $H$, a contradiction.

On the other hand, if $N_{F_{1}}(v_{T}) = N_{F_{1}}(v_{T}')$ for some $v_{T}$ and $v_{T}'$ in $T$, then by Claim~\ref{claim1} and Claim~\ref{claim4}, there exists one vertex in $F_{1}$ which is not adjacent to either of $v_{T}$ or $v_{T}'$. So $v_{T}$ is adjacent to $v_{T}'$ since $\alpha(G) \leq 2$.

So $N_{F_{1}}(v_{T}) = N_{F_{1}}(v_{T}')$ if and only if $v_{T}$ is adjacent to $v_{T}'$. Thus, $T$ is formed by some disjoint cliques. By $\alpha(G) \leq 2$, $T$ consists of at most two disjoint cliques, denoted by $T_{1}$ and $T_{2}$ ($T_{2}$ can be an empty set).

Now we get three cliques $F_{1}$, $T_{1} \cup F_{2}$ and $T_{2} \cup F_{2}$ in $G$, and $|F_{1}| + |T_{1} \cup F_{2}| + |T_{2} \cup F_{2}| = |G| + 3$. This indicates $\omega(G) \geq \lceil (|G| + 3) / 3 \rceil \geq \lceil (|G| + 3) / 4 \rceil$. By Theorem~\ref{packingseagull}, we get $h(G) \geq \chi(G)$, a contradiction.

\begin{case}
$|F_{2}| = 2$.
\end{case}

It is easy to see $|T| \geq 6$ since $|G| \geq 27$. Then $T$ contains a triangle $T_{1}$ since Ramsey number $R(3,3)=6$ and $\alpha(G) \leq 2$. By Claim~\ref{claim1} we get that for any vertex $v_{T} \in T$, $|N_{F_{1}}(v_{T})| \leq 3$. So $e(T_{1},F_{1}) \leq 3 \times 3 =9$. By pigeonhole principle, there exists a vertex $v_{F_{1}} \in F_{1}$, such that $|N_{T_{1}}(v_{F_{1}})| \leq 1$. If $|N_{T_{1}}(v_{F_{1}})| = 1$, then $\{ v_{F_{1}} \} \cup T_{1} \cup F_{2}$ form an induced copy of $H$, a contradiction. So $|N_{T_{1}}(v_{F_{1}})| = 0$.

Set $T_{2} = \{ v \in T \ |\ v \text{ is not adjacent to } v_{F_{1}} \}$, then $T_{1} \subseteq T_{2}$ and $T_{2}$ is a clique since $\alpha(G) \leq 2$. Set $T_{3} = T \backslash T_{2}$.

If $T_{3} = \emptyset$, then $F_{1}$ and $T \cup F_{2}$ are two disjoint cliques in $G$. Since $G = F_{1} \cup T \cup F_{2}$, $G$ contains a clique of at least $\lceil |G| / 2 \rceil$ vertices. Then $h(G) \geq \chi(G)$ by Theorem~\ref{HC2iff}, a contradiction.

If $T_{3} \neq \emptyset$. For any vertex $v_{T_{3}} \in T_{3}$ we have $|N_{T_{2}}(v_{T_{3}})| \leq 1$, otherwise $\{ v_{F_{1}} \} \cup \{ v_{T_{3}} \} \cup F_{2} \cup \{ \text{any two vertices in } N_{T_{2}}(v_{T_{3}}) \}$ form an induced copy of $H$. By $\alpha(G) \leq 2$ and $|T_{2}| \geq 3$, we get that $T_{3}$ is a clique. Now $|F_{1}| + |T_{2} \cup F_{2}| + |T_{3} \cup F_{2}| = |G| + 2$, we have $\omega(G) \geq \lceil (|G| + 2) / 3 \rceil \geq \lceil (|G| + 3) / 4 \rceil$. Then $h(G) \geq \chi(G)$ by Theorem~\ref{packingseagull}, a contradiction.

\begin{case}
$|F_{2}| = 1$.
\end{case}

It is easy to see $|T| \geq 9$ since $|G| \geq 27$. Then $T$ contains a $K_{4}$ since Ramsey number $R(3,4)=9$ and $\alpha(G) \leq 2$. If $|N_{K_{4}}(v_{F_{1}})| = 1$ for some vertex $v_{F_{1}} \in F_{1}$, then $\{ v_{F_{1}} \} \cup K_{4} \cup F_{2}$ form an induced copy of $H$. So for every vertex $v_{F_{1}} \in F_{1}$, $|N_{K_{4}}(v_{F_{1}})| \neq 1$.

\begin{subcase}
There exists a vertex $v_{F_{1}} \in F_{1}$, $|N_{K_{4}}(v_{F_{1}})| = 0$.
\end{subcase}

Set $T_{1} = \{ v \in T \ |\ v \text{ is adjacent to } v_{F_{1}} \}$, $T_{2} = T \backslash T_{1}$, then $T_{2}$ contains a $K_{4}$ (so $|T_{2}| \geq 4$) and $T_{2}$ is a clique since $\alpha(G) \leq 2$. If $|N_{T_{2}}(v_{T_{1}})| \geq 3$ for some vertex $v_{T_{1}} \in T_{1}$, then $\{ v_{T_{1}} \} \cup F_{2} \cup \{ v_{F_{1}} \} \cup \{ \text{any three vertices in } N_{T_{2}}(v_{T_{1}}) \}$ form an induced copy of $H$. So for every vertex $v_{T_{1}} \in T_{1}$, $|N_{T_{2}}(v_{T_{1}})| \leq 2$ .

\begin{subsubcase}
$|T_{2}| \geq 5$.
\end{subsubcase}

$|N_{T_{2}}(v_{T_{1}})| + |N_{T_{2}}(v_{T_{1}}')| \leq 4 < |T_{2}|$ for any two vertices $v_{T_{1}}$ and $v_{T_{1}}'$ in $T_{1}$. So there exists a vertex in $T_{2}$ that is not adjacent to either $v_{T_{1}}$ or $v_{T_{1}}'$. By $\alpha(G) \leq 2$, $v_{T_{1}}$ is adjacent to $v_{T_{1}}'$ hence $T_{1}$ is a clique due to the arbitrary choice of $v_{T_{1}}$ and $v_{T_{1}}'$. Now $|F_{1}| + |T_{1} \cup F_{2}| + |T_{2} \cup F_{2}| = |G| + 1$, we have $\omega(G) \geq \lceil (|G| + 1) / 3 \rceil \geq \lceil (|G| + 3) / 4 \rceil$. By Theorem~\ref{packingseagull} we get $h(G) \geq \chi(G)$, a contradiction.

\begin{subsubcase}
$|T_{2}| = 4$.
\end{subsubcase}

For any vertex $v_{T_{1}} \in T_{1}$, $N_{T_{2}}(v_{T_{1}}) \neq \emptyset$, otherwise $T_{2} \cup F_{2} \cup \{ v_{T_{1}} \}$ form an induced copy of $H$. Hence $1 \leq |N_{T_{2}}(v_{T_{1}})| \leq 2$ for every vertex $v_{T_{1}} \in T_{1}$.

If $T_{1}$ contains a $K_{5}$, then for every vertex $v_{T_{2}} \in T_{2}$, $|N_{K_{5}}(v_{T_{2}})| \neq 0,1,2$. Otherwise $\{ v_{T_{2}} \} \cup K_{5} \cup F_{2} \cup \{ v_{F_{1}} \}$ contain an induced copy of $H$. So $e(T_{2},K_{5}) \geq 3 \times 4 =12$. However, $|N_{T_{2}}(v_{T_{1}})| \leq 2$ for every vertex $v_{T_{1}} \in T_{1}$, we have $e(T_{2},K_{5}) \leq 2 \times 5 =10$, a contradiction.

If $T_{1}$ is $K_{5}$-free, by Ramsey number $R(3,5) = 14$ and $\alpha(G) \leq 2$, we obtain that $|T_{1}| \leq 13$. So $|G| \leq 24$, a contradiction.

\begin{subcase}
For every vertex $v_{F_{1}} \in F_{1}$, $|N_{K_{4}}(v_{F_{1}})| \geq 2$.
\end{subcase}

\begin{subsubcase}\label{case3.2.1}
$|G| \geq 29$
\end{subsubcase}

$|T \backslash K_{4}| \geq 18$ since $|G| \geq 29$. By Ramsey number $R(3,6) = 18$ and $\alpha(G) \leq 2$, $T \backslash K_{4}$ contains a $K_{6}$.

For every vertex $v_{K_{4}} \in K_{4}$, $|N_{K_{6}}(v_{K_{4}})| \neq 0,1,2$, otherwise $\{ v_{K_{4}} \} \cup K_{6} \cup F_{2}$ contain an induced copy of $H$. For every vertex $v_{F_{1}} \in F_{1}$, $|N_{K_{6}}(v_{F_{1}})| \neq 1,2,3$, otherwise $\{ v_{F_{1}} \} \cup K_{6} \cup F_{2}$ contain an induced copy of $H$.

If there exists a vertex $v_{F_{1}} \in F_{1}$ with $|N_{K_{6}}(v_{F_{1}})| = 0$, then we can choose a vertex $v_{K_{4}} \in N_{K_{4}}(v_{F_{1}})$ such that $\{ \text{any three vertices in } N_{K_{6}}(v_{K_{4}}) \} \cup \{ v_{F_{1}} \} \cup \{ v_{K_{4}} \} \cup F_{2}$ forms an induced copy of $H$.

If for every vertex $v_{F_{1}} \in F_{1}$, $|N_{K_{6}}(v_{F_{1}})| \geq 4$, then $e(F_{1},K_{6}) \geq 4 \times 6 = 24$. However, $|N_{F_{1}}(v_{T})| \leq 3$ for every vertex $v_{T} \in T$ by Claim~\ref{claim1}, we have $e(F_{1},K_{6}) \leq 3 \times 6 = 18$, a contradiction.

\begin{subsubcase}
$|G| = 27$
\end{subsubcase}

$|T| = 20$ since $|G| = 27$. Let $A$ be a set of two vertices in $K_{4}$, by Ramsey number $R(3,6) = 18$ and $\alpha(G) \leq 2$, $T \backslash A$ contains a $K_{6}$. We have $|K_{6} \cap K_{4}| \leq 2$.

If $|K_{6} \cap K_{4}| = 0$, the remaining proof is the same as Case~\ref{case3.2.1}.

If $1 \leq |K_{6} \cap K_{4}| \leq 2$, then we can divide the proof into two parts:

If there exists a vertex $v_{F_{1}} \in F_{1}$ with $N_{K_{6}}(v_{F_{1}}) = \emptyset$, then we can choose a vertex $v_{K_{4}} \in N_{K_{4}}(v_{F_{1}})$ such that $\{ v_{F_{1}} \} \cup \{ v_{K_{4}} \} \cup K_{6} \cup F_{2}$ contain an induced copy of $H$ (either $v_{F_{1}}$ or $v_{K_{4}}$ is the vertex of degree one in $H$).

If $N_{K_{6}}(v_{F_{1}}) \neq \emptyset$ for every vertex $v_{F_{1}} \in F_{1}$, then we can also get that $|N_{K_{6}}(v_{F_{1}})| \neq 1,2,3$ for every vertex $v_{F_{1}} \in F_{1}$, otherwise $\{ v_{F_{1}} \} \cup K_{6} \cup F_{2}$ contain an induced copy of $H$. So $e(F_{1},K_{6}) \geq 4 \times 6 = 24$. On the other hand, $|N_{F_{1}}(v_{T})| \leq 3$ for every vertex $v_{T} \in T$ by Claim~\ref{claim1}, we have $e(F_{1},K_{6}) \leq 3 \times 6 = 18$, a contradiction.\\

So far we have proved the Theorem~\ref{T2} for $H = K_{5}^{+}$, it suffices to prove that Theorem~\ref{T2} is also holds for $H = H_{7}$.

Let $H = H_{7}$ and let $G$ be a minimal counterexample to Theorem~\ref{T2} with respect to $|V(G)| + |E(G)|$. That is to say, $G$ is $H$-free with $\alpha(G) \leq 2$, but $h(G) < \chi(G)$.

By the proof above, $G$ contains an induced copy of $K_{5}^{+}$. Let $v$ be the vertex of degree one in $K_{5}^{+}$, and $u$ be the neighbor of $v$ in $K_{5}^{+}$.

For every vertex $w \in G \backslash K_{5}^{+}$, if $w$ is not adjacent to $v$, then $w$ is adjacent to all vertices of $K_{5}^{+} \backslash \{ u,v \}$ since $\alpha(G) \leq 2$. So $w$ is adjacent to $u$ because $G$ is $H$-free. That is to say, edge $e_{uv}$ is a dominating edge, a contradiction to Lemma~\ref{lemma2}.

This completes the proof of Theorem~\ref{T2}.\\

\section{Proof of Theorem~\ref{T3}}

By Ramsey number $R(3,8) = 28$ and $\alpha(G) \leq 2$, every $K_{8}$-free counterexample to Problem~\ref{HC2} has exactly 27 vertices. By Ramsey number $R(3,7) = 23$ and $K_{8}$-free we have $\Delta(G) \leq 22$ and $\delta(G) \geq 19$, otherwise the neighbors of a vertex $v$ with maximum degree has a $K_{7}$ subgraph(together with $v$ we get $K_{8}$), or the non-neighbors of a vertex $v$ with minimum degree has a $K_{8}$ subgraph.

Let $\Delta(G) = 22$ and $v$ be the vertex of degree 22. Let $A = V(G)\backslash N_{G}[v]$, then $|A| = 4$. For any vertex $u \in N_{G}(v)$, $u$ is not adjacent to at least one vertex in $A$, otherwise $e_{uv}$ is a dominating edge. Set $A = \{ v_{1}, v_{2}, v_{3}, v_{4} \}$ and $N_{G}(v) = N_{1} \cup N_{2} \cup N_{3} \cup N_{4}$ where any vertex in $N_{i}$ is adjacent to $v_{1},...,v_{i-1}$ but not adjacent to $v_{i}$ for $i = 1,2,3,4$. Then $N_{i}$ is a clique since $\alpha(G) \leq 2$, so $|N_{1}| \leq 6, |N_{2}| \leq 6, |N_{3}| \leq 5, |N_{4}| \leq 4$ since $G$ is $K_{8}$-free. But $22 = |N_{G}(v)| = |N_{1}| + |N_{2}| + |N_{3}| + |N_{4}| \leq 21$, a contradiction.

By the assumption $\delta(G) \neq 19$ we have $20 = \delta(G) \leq \Delta(G) \leq 21$ since there is no 21 regular graph with 27 vertices.

Let $\Delta(G) = 21$ and $v$ be the vertex of degree 21. Let $A = V(G)\backslash N_{G}[v]$, then $|A| = 5$. For any vertex $u \in N_{G}(v)$, $u$ is not adjacent to at least one vertex in $A$, otherwise $e_{uv}$ is a dominating edge. So $80 = 5 \times (20 - 4) \leq e(N_{G}(v),A) \leq 21 \times 4 =84$. Set $A = \{ v_{1}, v_{2}, v_{3}, v_{4}, v_{5} \}$ and $N_{G}(v) = N_{1} \cup N_{2} \cup N_{3} \cup N_{4} \cup N_{5} \cup N$, where any vertex in $N_{i}$ is not adjacent to $v_{i}$ and adjacent to $v_{j}$($j \neq i$) for $i = 1,2,3,4,5$, any vertex in $N$ is not adjacent to at least two vertices in $A$. Then $N_{i}$ is a clique since $\alpha(G) \leq 2$, so $|N_{i}| \leq 3$ since $G$ is $K_{8}$-free. But $|N| \leq 4$ since $e(N_{G}(v),A) \geq 80$, we have $22 = |N_{G}(v)| = |N_{1}| + |N_{2}| + |N_{3}| + |N_{4}| + |N_{5}| +|N| \leq 19$, a contradiction.

If $\Delta(G) = 20$, let $v$ be a vertex of degree 20. Let $A = V(G)\backslash N_{G}[v]$, then $|A| = 6$. For any vertex $u \in N_{G}(v)$, $u$ is not adjacent to at least one vertex in $A$, otherwise $e_{uv}$ is a dominating edge. Set $N_{G}(v) = N_{1} \cup N_{2}$, where any vertex in $N_{1}$ is adjacent to five vertices in $A$ and any vertex in $N_{2}$ is adjacent to at most four vertices in $A$. Since $e(N_{G}(v),A) = 6 \times (20 - 5) = 90$, we have $|N_{2}|\leq 10$. There are at most two vertices in $N_{1}$ have common neighbors in $A$ since $G$ is $K_{8}$-free, so $|N_{1}| \leq 12$. Thus $N_{2}$ is one of the following:
\begin{itemize}
\item[$\bullet$] 10 vertices which have 4 neighbors in $A$; or
\item[$\bullet$] 9 vertices which have 4 neighbors in $A$; or
\item[$\bullet$] 8 vertices which have 4 neighbors in $A$; or
\item[$\bullet$] 8 vertices which have 4 neighbors in $A$ and 1 vertex which has 3 neighbors in $A$; or
\item[$\bullet$] 7 vertices which have 4 neighbors in $A$ and 1 vertex which has 3 neighbors in $A$; or
\item[$\bullet$] 7 vertices which have 4 neighbors in $A$ and 1 vertex which has 2 neighbors in $A$; or
\item[$\bullet$] 6 vertices which have 4 neighbors in $A$ and 2 vertices which have 3 neighbors in $A$.
\end{itemize}

The corresponding number of edges in $A$ that cannot dominate $N_{2}$ is at most 10, 9, 8, 11, 10, 13, 12. So there exists one edge $e = \{ u,w\}$ in $A$ which dominates $N_{2}$ since $e(A) = 15$. By the definition of $N_{1}$, $e$ dominates $N_{1}$ as well.

The vertex $u$ or $w$ has at least one non-neighbor in $N_{1}$, otherwise there are 3 vertices in $N_{1}$ have common neighbors in $A$ since $|N_{1}| \geq 10$, a contradiction to $K_{8}$-free. Without loss of generality, we assume that $x \in N_{1}$ is not adjacent to $u$, then $e' = \{ v,x\}$ and $e = \{ u,w\}$ are connected dominating matching, a contradiction to~\ref{lemma2}.

This completes the proof of Theorem~\ref{T3}.\\

\bibliographystyle{abbrv}
\bibliography{mybib}

\end{document}